\documentclass{amsart}
\usepackage{graphicx}
\usepackage{geometry}

\usepackage{amssymb}
\usepackage{amsmath}
\usepackage{array}
\usepackage{enumerate}
\usepackage{multirow}
\usepackage{float}
\usepackage{hyperref}
\usepackage[main=british, french, german, latin]{babel}
\newcommand{\Hy}{\mathbb{H}}
\newcommand{\Euc}{\mathbb{R}}

\newcommand{\Sph}{\mathbb{S}}

\newcommand{\lla}{\left\langle}
\newcommand{\rra}{\right\rangle}

\newcommand{\Pf}{{\noindent\em \textbf{Proof}}. }
\newcommand{\PfRef}[1]{{\noindent\em \textbf{Proof of Theorem #1}}. }
\newcommand{\EPf}{\hfill$\square$\newline}

\newcommand{\surface}[1]{\left(\Sigma^2, #1 \right)}

\newcounter{theorem}

\newtheorem{thm}[theorem]{Theorem}
\newtheorem*{thmRef}{Theorem}

\newtheorem{cor}[theorem]{Corollary}

\newtheorem{lem}[theorem]{Lemma}

\newtheorem{dfn}[theorem]{Definition}

\theoremstyle{remark}

\title{A note on parallel mean curvature surfaces and Codazzi operators}

\author{Felippe Guimar\~aes}%\thanks{}
\address{KU Leuven, Department of Mathematics, Celestijnenlaan 200 B – Box 2400, 3001 Leuven, Belgium.}
\email{felippe.guima@gmail.com}

\subjclass[2010]{53C20 (Primary); 53A10, 53C42 (Secondary)}

\begin{document}

	\maketitle
	
    \begin{abstract}
        We present an intrinsic Klotz-Osserman type theorem for surfaces in terms of Codazzi operators. Additionally, utilizing Simons' formula, we investigate surfaces with parallel mean curvature with non-positive Gaussian curvature in product spaces.
    \end{abstract}

    \section{Introduction}

    Consider $\Sigma^2$ as an orientable $2$-dimensional Riemannian manifold. We define a symmetric endomorphism $S$ of the tangent space $T\Sigma$ to be a Codazzi operator if it satisfies the Codazzi equation given by $$ \nabla_X SY - \nabla_Y SX - S[X,Y] =0\ \text{for all}\ X,Y \in T\Sigma,$$ where $\nabla$ is the Levi-Civita connection of the metric of $\Sigma^2$. In the context of the shape operator for surfaces in three-dimensional space forms, such as $\Euc^3$, $\Sph^3$, and $\Hy^3$,  the Codazzi equation serves as one of the two primary integrability conditions, alongside the Gauss equation. Some classical results in surface theory in $\Euc^3$ depend essentially on the Codazzi equation alone. Among those are the classical Hopf's theorem, which asserts the uniqueness of round spheres among all immersed spheres with constant mean curvature. There's also Efimov's theorem, which states that no complete surface immersed in $\Euc^3$ can have a negative Gaussian curvature that stays away from zero (see \cite{SmythXavierEfimovs}, \cite{TillaEfimov}, \cite{Efimov}). This led to the adaptation of numerous results into an abstract setting of Codazzi operators, a topic well explored in \cite{AledoEspinarCodazziPairs}. Several elegant proofs emerged from this intrinsic approach in works such as \cite{GalvezMartinezTeruelSurfaceCodazzi}, where the authors explore the existence of various Codazzi operators on hypersurfaces in space forms, including the second fundamental form of these hypersurfaces and the equidistant hypersurfaces. However, in more general settings, the second fundamental form may not be a Codazzi operator.

    Nonetheless, when studying surfaces with constant mean curvature (CMC surfaces) in different ambient spaces, there exists a vast literature on holomorphic quadratic differentials after Hopf discovered that the complexification of the traceless part of the second fundamental form and the seminal work \cite{AbreschRosenbergHolomorphicForm} (see also\cite{EspinarTrejoTeseHaimer, FectuPmcSurfacesComplexSpaceForms} and references therein). These forms give rise to an endomorphism, which becomes a natural candidate for a Codazzi operator, as indicated in \cite{AledoEspinarCodazziPairs}.
    
    Among the classical results on CMC surfaces, we highlight the Klotz-Osserman's theorem.
    
    \begin{thmRef}[\cite{KlotzOsserman}]
    A complete CMC surface in $\Euc^3$ whose Gaussian curvature $K$ does not change sign is either a sphere, a minimal surface, or a right circular cylinder.
    \end{thmRef}
    
    This result was extended to $\Sph^3$ in \cite{HoffmanCMCSphere} and to $\Hy^3$ in \cite{TribuzyHypHopf}, with an additional hypothesis when the curvature is non-positive. Over a decade ago, in \cite{EspinarRosenbergCMCHomogeneous}, the authors obtained a version of this result for the homogeneous three-manifolds with a four-dimensional isometry group.
    
    Inspired by the works \cite{EspinarRosenbergCMCHomogeneous}, \cite{EspinarTrejoTeseHaimer} and \cite{GalvezMartinezTeruelSurfaceCodazzi}, we formulate a Klotz-Osserman type theorem for Codazzi operators.

    \begin{thm}\label{thm:MainCodazzi}
    Let $\Sigma^2$ be an orientable surface with Gaussian curvature $K$ that is not topologically a sphere with a traceless Codazzi operator $S$. Then, 
    \begin{itemize}
        \item\label{thmI:pos} If $K \geq 0$ and $|\det S| \leq \epsilon$ for some constant $\epsilon$, then $\det S$ is constant and $K.\det S \equiv 0$.
        \item\label{thmI:Kneg} If $K \leq 0$ and $|\det S| \geq \epsilon > 0$ for some constant $\epsilon$, then $\det S$ is constant and $K \equiv 0$.
    \end{itemize}
    \end{thm}

    As an application of Theorem \ref{thm:MainCodazzi}, we investigate surfaces with parallel mean curvature (PMC surfaces) that have non-positive Gaussian curvature in the product of a space form $\mathbb{M}^n(\kappa)$ with a line. The objective is to address the case not covered by \cite{AlencarTribuzyHopfMxR} and \cite{FectuRosenbergNotePMCsurfacesProducts}, specifically when the PMC surface exhibits non-positive Gaussian curvature.
    
    \begin{thm}\label{thm:PMCSurfacesProductsMu}
    Let $\Sigma^2 \subset \mathbb{M}^n(\kappa)\times \Euc$, where $\kappa \neq 0$, be a complete non-minimal PMC surface with Gaussian curvature $K \leq 0$ such that  $\mu = \sup_{\Sigma}(|\alpha|^2 - \frac{1}{|H|^2}|A_H|^2) < \infty$.  Then, $K \equiv 0$ if
    \begin{itemize}
        \item $\kappa <0$:\quad $K - \vert H \vert^2 -\kappa + \frac{1}{2}\mu \leq -\epsilon <0$ for some constant $\epsilon>0$;
        \item $\kappa >0$:\quad$K + \frac{1}{2}\mu \leq c \vert H \vert^2$ and $4(1-c)^{\frac{1}{2}}\vert H \vert^2 -\kappa>0$ for some constant $0\leq c <1$.
    \end{itemize}
    \end{thm}

    The constant $\mu$ naturally emerges with the increase of the codimension, reflecting the need for curvature control of the normal bundle (see \cite{FectuRosenbergPMCemQ3R}, \cite{BatistaCavalcanteFectuFiniteTotalCurvatureProducts}). This constant is related to how close the surface is to being umbilical. For minimal surfaces, considering only the angle that the surface makes with the parallel field of the ambient space, we can construct a Codazzi operator. In this scenario, if the surface is far from being a slice, then the angle function must be constant (see Corollary \ref{cor:AngleProducts}).

    The structure of the article is outlined as follows: In Section \ref{sec:Preliminaries}, we provide the necessary definitions and tools used throughout this manuscript. In Section \ref{sec:Results}, the proofs of the main theorems are presented, along with applications for surfaces in product spaces. Lastly, in Section \ref{sec:Remarks}, we make some remarks about the assumptions and results obtained.

    \section{Preliminaries}\label{sec:Preliminaries}
    Let $\Sigma^2 \subset N^n\times \Euc$ be a complete orientable surface immersed in the product of an $n$-dimensional Riemannian manifold $N^n$ and $\Euc$. This surface inherits an induced metric $\lla \cdot, \cdot \rra$, and for simplicity, we will use the same notation to denote both the metric of the ambient space and that of the surface. We denote the second fundamental form of $\Sigma^2$ by $\alpha(X,Y) = ({}^N\nabla_X Y)|_{T\Sigma}$, where ${}^N\nabla_X Y$ is the Levi-Civita connection of the ambient space. For this work, we will make use of the natural orthogonal decomposition of the tangent space $T(N^n\times \Euc) = T\Sigma \oplus T^{\perp}\Sigma$ over $\Sigma^2$. Let $\xi \in T^{\perp}\Sigma$ be a normal vector field. The Weingarten equation is given by $$\nabla_X \xi = -A_{\xi} X + \nabla^{\perp}_X \xi,$$ where $X \in T\Sigma$, $\nabla^{\perp}$ is the normal connection, and $A_{\xi}$ is the shape operator associated with $\xi$, that is, $\lla A_{\xi} X, Y\rra = \lla \alpha(X,Y), \xi \rra$ for all $X,Y \in T\Sigma$. Given an orthonormal basis ${X,Y}$ of $T_p\Sigma$ at point $p \in \Sigma^2$, the mean curvature vector of the surface $\Sigma^2$ is defined as $H(p) = \frac{1}{2} (\alpha(X,X) + \alpha(Y,Y))(p)$.

    \begin{dfn}
        The surface $\Sigma^2\subset N^n\times \Euc$ is called a parallel mean curvature surface, or just PMC surface, if its mean curvature vector $H$ is parallel in the normal bundle, i.e., $\nabla^{\perp} H = 0$.
    \end{dfn}

    Consider the vector field  $\frac{\partial}{\partial t}$ that generates the line $\Euc$ of $N^n\times \Euc$. Decomposing this with respect to the surface $\Sigma^2$, we obtain the orthogonal decomposition $$\frac{\partial}{\partial t} = T + \eta,$$ where $T \in T\Sigma^2$ and $\eta \in T^{\perp}\Sigma$. Due to the parallelism of $\frac{\partial}{\partial t}$, we derive the relations $\nabla_X T = A_{\eta}T$ and $\alpha(X,T) = -\nabla^{\perp}_X \eta$. When $N^n = \mathbb{M}^n(\kappa)$ is a space form (a simply connected manifold with constant sectional curvature $\kappa$), we have the following fundamental equations.
    
    \noindent The Gauss equation $$R(X,Y)W = \epsilon (X \wedge Y - \lla Y, T\rra X \wedge T + \lla X,T \rra Y \wedge T)W + A_{\alpha(Y,W)} X - A_{\alpha(X,W)}Y$$ and the Codazzi equation 
    \begin{equation}\label{eq:CodazziProduct}
        (\nabla_Y A)(X,\xi) - (\nabla_X A)(Y,\xi) = \kappa \lla \xi, \eta \rra (X \wedge Y)T,
    \end{equation} where $(X \wedge Y)T = \lla Y, T \rra X - \lla X, T\rra Y$ and $X,Y \in T\Sigma$.

    The work of \cite{AlencarTribuzyHopfMxR} introduced a real quadratic form on PMC surfaces, leading to the traceless endomorphism
    
    \begin{equation}\label{eq:CodazziPMCMxR}
        SX = 2A_{H} -\kappa\lla T,X \rra T + \kappa \frac{|T|^2}{2}X -2|H|^2X.
    \end{equation}

    This operator, first introduced in \cite{BatistaCodazziOperatorProducts} and later used by \cite{FectuRosenbergNotePMCsurfacesProducts}, has been instrumental in the study of PMC surfaces with non-negative Gaussian curvature. Importantly, \cite{BatistaCodazziOperatorProducts} proved that $S$ is a Codazzi operator when the ambient space is $3$-dimensional. The same proof holds when the surface has arbitrary codimension, as detailed below.

    \begin{lem}\label{lem:SCodazzi}
        The operator $S$ given by \eqref{eq:CodazziPMCMxR} satisfies the Codazzi equation.
    \end{lem}
    \Pf
        Consider the operator $\tilde{T}_S(X,Y) = \nabla_X SY - \nabla_Y SX - S[X,Y]$ for $X,Y \in T\Sigma$. As $\tilde{T}_S$ is tensorial, bilinear, and skew-symmetric, it suffices to show that for every $p \in T\Sigma$, there exists a basis $\{X,Y\}$ of $T_p\Sigma$ such that $\tilde{T}_S(X,Y) = 0$.
    
        For $p \in \Sigma^2$, if $T = 0$, then it directly follows from \eqref{eq:CodazziProduct} that $\tilde{T}_S = 0$. Now suppose the contrary. Let $\{X,Y\}$ be an orthonormal basis that diagonalizes $A_{\eta}$ with eigenvalues $\lambda_X$ and $\lambda_Y$. Then we have the following expression:
    
        \begin{equation}
            \begin{split}
                \tilde{T}_S(X,Y) & = \kappa \left[ 2 \lla H, \eta\rra (Y \wedge X) T -\lla \nabla_X T, Y \rra T - \lla T, Y \rra \nabla_X T + \lla \nabla_X T, T\rra Y \right. \\
                & \quad \left. +\lla \nabla_Y T, X \rra T + \lla T, X \rra \nabla_Y T - \lla \nabla_Y T, T\rra X \right] \\
                & = \kappa \left[ 2 \lla H, \eta\rra (Y \wedge X) T - (\lambda_X + \lambda_Y) (Y \wedge X) T \right] \\
                & = 0,
            \end{split}
        \end{equation}
        where in the last equality we used that $(\text{trace~}A_\eta) = 2\lla H, \eta \rra$. Thus, the lemma follows.
    \EPf

    Consider a local orthonormal frame field in the normal bundle given by $\{\xi_1 = \frac{H}{|H|}, \xi_2, \dots, \xi_{n-1}\}$ and denote $A_i = A_{\xi_i}$. Note that trace $A_1=2|H|$ and trace $A_i=0$ for all $i>1$. Following a straightforward computation (see \cite{BatistaCavalcanteFectuFiniteTotalCurvatureProducts}, \cite{FectuRosenbergNotePMCsurfacesProducts}), we derive the expression for the non-minimal PMC surface Gaussian curvature as

    \begin{equation}\label{eq:ProductKg}
        K = \kappa(1-| T |^2) + |H|^2 - \frac{1}{8|H|^2} |S|^2 - \frac{\kappa^2}{16|H|^2}|T|^4 - \frac{\kappa}{4|H|^2}\lla ST, T \rra + \sum_{i>1} \det A_i.
    \end{equation}

    \noindent  For the proof of Theorem \ref{thm:MainCodazzi}, we will lean on the Simons's formula for Codazzi operators, stated below.

    \begin{thmRef}[Theorem 3 in \cite{EspinarTrejoTeseHaimer}, \cite{ChengYauCodazziSimon}]
        Let $\Sigma^2$ be an orientable surface equipped with a traceless Codazzi operator $S$. Then we have that $|S|^2 = -2\det S$ satisfies the equation
        \begin{equation}\label{eq:SimonsType}
            \frac{1}{2} \Delta|S|^2 = |\nabla |S||^2 + 2K |S|^2,
        \end{equation}
        which, away from zeros of $|S|$, can also be represented as
        \begin{equation}\label{eq:ReducedEquation}
            |S|\Delta|S| - 2K |S|^2 = |\nabla |S||^2,
        \end{equation}

        \noindent or equivalently,
        
        \begin{equation}\label{eq:LnS}
            \Delta\ln|S| = 2K.
        \end{equation}

        Here, $\Delta$ represents the Laplace–Beltrami operator relative to the metric of $\Sigma^2$.
    \end{thmRef}

     The subsequent lemma will be useful for the metric change argument. It arises from straightforward computations involving Koszul's formula and the use of the curvature tensor definition. The metric's completeness is deduced using the well-known Hopf-Rinow theorem.

    \begin{lem}\label{lem:ChangeMetrics}
    Let $\Sigma^2$ be a complete Riemannian manifold endowed with the metric $\lla \cdot, \cdot \rra$, the Levi-Civita connection $\nabla$, and Gaussian curvature $K$. Suppose there exists a traceless Codazzi operator $S$ such that $|S| \geq \epsilon>0$ for some constant $\epsilon$. Under these conditions, the newly defined metric $\lla \cdot, \cdot\rra_S = \lla S \cdot, S \cdot \rra$ is complete. Furthermore, the Levi-Civita connection $\Tilde{\nabla}_X Y$ corresponding to this new metric is $$\Tilde{\nabla}_X Y = S^{-1}\nabla_X SY.$$ In particular, the Gaussian curvature $\Tilde{K}$ of $\Sigma^2$ associated with the metric $\lla \cdot, \cdot \rra$ is given by $\Tilde{K} = \frac{K}{\det S}$.
    \end{lem}

    \section{The results}\label{sec:Results}

    In this section, we present the proof of the main theorem and discuss some applications. We utilize Simons' formula, as introduced in Section \ref{sec:Preliminaries}, alongside techniques from \cite{HuberSubharmonic,KlotzOsserman}. Our strategy aligns with the traditional methods seen in Osserman-Klotz-type theorems, as in \cite{EspinarRosenbergCMCHomogeneous,FectuRosenbergNotePMCsurfacesProducts}.

    \PfRef{\ref{thm:MainCodazzi}}
    Assume that we are in the situation that $K \geq 0$. By Lemma 5 in \cite{KlotzOsserman}, it follows that $\Sigma^2$ is either a sphere or non-compact and parabolic. In the latter scenario, from equation \eqref{eq:SimonsType}, we deduce that $|S|^2 = - 2\det S$ is a subharmonic function bounded from above on a non-compact parabolic surface. Hence, by \cite[page~204]{AhlforsSarioBookRiemannSurfaces} (or see Proposition 3 and Corollary 1 on page~336 in \cite{ChengYauSubharmonic}), $\det S$ must be constant and \eqref{eq:SimonsType} implies that $K.\det S \equiv 0$.
    
    Moving on to the case where $K \leq 0$, since $S$ is non-singular, we can modify the surface metric to $\surface{\lla\cdot, \cdot\rra_S}$, where $\lla X,Y \rra_S = \lla SX,SY \rra$ for $X,Y \in T\Sigma$. Lemma \ref{lem:ChangeMetrics} ensures that this new metric is complete and satisfies $\Tilde{K} \geq 0$. Using Lemma 5 in \cite{KlotzOsserman} again, we deduce that $\surface{\lla\cdot, \cdot\rra_S}$ is parabolic. Considering the inequality
    $$\Tilde{\Delta} \ln |S| = \frac{4K}{|S|^2} \leq 0,$$
    where $\Tilde{\Delta}$ is the Laplacian associated with the metric $\lla \cdot, \cdot\rra_S$, and noting that $\ln |S|$ has a lower bounded, we conclude that both $\det S$ and $K$ are constants.
    \EPf

    Consider a surface immersed in an arbitrary Riemannian manifold. One method to construct a Codazzi operator involves using the fundamental equations of submanifolds. However, a clever way to find a candidate is through holomorphic quadratic forms. The endomorphism associated with this quadratic form emerges as a suitable choice, as evidenced by Lemma 1 in \cite{AledoEspinarCodazziPairs}. This approach is employed in \cite{EspinarTrejoTeseHaimer, FectuRosenbergNotePMCsurfacesProducts}. To apply Theorem \ref{thm:MainCodazzi}, we will impose curvature constraints and deduce estimates for $|S|^2$ using \eqref{eq:ProductKg}.

    The first application concerns PMC surfaces with non-positive curvature in the ambient space $\mathbb{M}^n(\kappa) \times \Euc$. While \cite{AlencarTribuzyHopfMxR} addressed the case where $K \geq 0$, the scenario with non-positive sectional curvature had not been previously explored. We present this result as an application of Theorem \ref{thm:MainCodazzi}. In particular, it implies Theorem~\ref{thm:PMCSurfacesProductsMu}.

    \begin{thm}\label{thm:PMCSurfacesProducts}
        Let $\Sigma^2 \subset \mathbb{M}^n(\kappa)\times \Euc$, where $\kappa \neq 0$, be an orientable complete non-minimal PMC surface with $K \leq 0$. Then, $K \equiv 0$ if 
        \begin{itemize}
            \item $\kappa<0$: $K - \vert H \vert^2 -\kappa - \sum_{i>1} \det A_{i} \leq -\epsilon <0$ for some constant $\epsilon>0$,
            \item $\kappa>0$: $K - \sum_{i>1} \det A_{i} \leq c \vert H \vert^2$ and $4(1-c)^{\frac{1}{2}}\vert H \vert^2 -\kappa>0$ for some constant $0\leq c <1$.
        \end{itemize}
    \end{thm}
    
    \Pf
        The objective is to show that $|S|$ is bounded from below and then apply Theorem \ref{thm:MainCodazzi}.
    
        For the case when $\kappa<0$: From \eqref{eq:ProductKg}, we have 
        \begin{equation}
        \frac{1}{8|H|^2} |S|^2 + \frac{|\kappa|}{4|H|^2} |S| \geq - \left( \frac{\kappa}{4|H|}|T|^2 + 2|H| \right)^2 +4|H|^2 - (K - \vert H \vert^2 -\kappa - \sum_{i>1} \det A_{i}).
        \end{equation}
        Using the given hypothesis, we obtain 
        \begin{equation}
        \frac{1}{8|H|^2} |S|^2 + \frac{|\kappa|}{4|H|^2} |S| \geq \epsilon > 0,
        \end{equation}
        which implies that $|S|$ cannot be arbitrarily small.
    
        For the case when $\kappa>0$. We no longer have the inequality $- \left( \frac{\kappa}{4|H|}|T|^2 + 2|H| \right)^2 +4|H|^2 > 0$. However, from \eqref{eq:ProductKg}, we deduce 
        \begin{equation}
        \frac{1}{8|H|^2} |S|^2 + \frac{|\kappa|}{4|H|^2} |S| +(K - \sum_{i>1} \det A_{i}) \geq \kappa (1 -|T|^2) + |H|^2 - \frac{\kappa^2}{16|H|^2} |T|^4.
        \end{equation}
        Utilizing the condition $K - \sum_{i>1} \det A_{i} \leq c \vert H \vert^2$, we get 
        \begin{equation}
        \frac{1}{8|H|^2} |S|^2 + \frac{|\kappa|}{4|H|^2} |S| \geq (1-c)|H|^2 - \frac{\kappa^2}{16|H|^2} = \frac{1}{16|H|^2}(4(1-c)^{\frac{1}{2}}|H|^2 - \kappa)(4(1-c)^{\frac{1}{2}}|H|^2 + \kappa),
        \end{equation}
        which, combined with the inequality $4(1-c)^{\frac{1}{2}}\vert H \vert^2 -\kappa>0$, concludes the proof.
    \EPf

    Note that in the case of hypersurfaces, we have $\mu = 0$, and the hypotheses in the case $\kappa <0$ are the same as those in \cite{TribuzyHypHopf}. Meanwhile, the hypotheses for the situation where $\kappa>0$ in the aforementioned theorem draw inspiration from those outlined in Remark 4.1 of \cite{EspinarRosenbergCMCHomogeneous} when the ambient space is a product manifold.

   Direct applications of Theorem \ref{thm:MainCodazzi} can be stated as follows.

    \begin{cor}\label{cor:MinimalSurfacesProducts}
        Let $\Sigma^2 \subset \mathbb{M}^n(\kappa)\times \Euc$, where $\kappa \neq 0$, be an orientable complete minimal surface. If $K \leq 0$ and $|T| > \epsilon$ for some constant $\epsilon>0$, then we have $K \equiv 0$ and $|T|$ is constant. Moreover, one of the following assertions holds:
            \begin{enumerate}
            \item\label{thmProducts:ItemVertical} $\Sigma^2$ is a vertical cylinder over a complete geodesic;
            \item\label{thmProducts:ItemExample} $\Sigma^2$ lies in $\mathbb{M}^4(\kappa)\times \Euc \subset \Euc^6$ with $\kappa >0$. Up to an isometry of $\mathbb{M}^4(\kappa)\times \Euc$, $\Sigma^2$ is parameterized by
                $$f(u,v) = \left( \frac{\cos \theta}{b}\cos{bu}, \frac{\cos \theta}{b}\sin{bu}, \frac{\sin{bv}}{b}, \frac{\cos{bv}}{b}, 0, u\sin \theta \right),$$ where $\sin{\theta} = |T|$ and $b = \sqrt{\kappa +\kappa \cos^2\theta}$.    
        \end{enumerate}
    \end{cor}
    \Pf
    It is straightforward to check that $|S|^2 = \frac{\kappa^2}{2}|T|^4$ and we can apply Theorem \ref{thm:MainCodazzi} and it follows that $K \equiv 0$ and $|T|$ is constant. The classification follows from the main theorem in \cite{HouQiuFlatProductConstantAngle}.
    \EPf
   
   \begin{cor}\label{cor:trivialCor}
        Let $\Sigma^2 \subset \mathbb{H}^n\times \Euc$, be an orientable complete minimal surface.  If $|T| > \epsilon$ for some constant $\epsilon>0$, then $\Sigma^2$ is a vertical cylinder over a complete geodesic.
    \end{cor}

    The same idea used in Corollary \ref{cor:MinimalSurfacesProducts} can be applied to a slightly more general ambient space. Let $\Sigma^2$ be a minimal surface in the product $N^n\times \Euc$. Just as in $\mathbb{M}^n(\kappa)\times\Euc$, we can decompose $\frac{\partial}{\partial t} = T+\eta$. Given that $\frac{\partial}{\partial t}$ is a parallel vector field, we obtain the same expressions $\nabla_X T = A_\eta X$ and $\alpha(X,T) = -\nabla^{\perp}_X \eta$ for all $X \in TM$. 
    
    We can then define the operator $\Tilde{S}X = -\lla T,X \rra T + \frac{|T|^2}{2} X$. Similar to the calculation in Lemma \ref{lem:SCodazzi}, it can be shown directly that \(\Tilde{S}\) is a Codazzi tensor. Building on Theorem \ref{thm:MainCodazzi}, we present the following application.

    \begin{cor}\label{cor:AngleProducts}
        Let \(\Sigma^2 \subset N^n\times \Euc\) be an orientable complete minimal surface. Then,
    
        \begin{itemize}
            \item If \(K \leq 0\) and \(|T| > \epsilon\) for some constant $\epsilon>0$ then \(K \equiv 0\) and \(|T|\) is constant. 
            \item If \(K \geq 0\) then either \(\Sigma^2\) is a topological sphere or \(|T|\) is constant and \(K.|T| \equiv 0\).
        \end{itemize}
    \end{cor}
    
    An important observation from the aforementioned corollary is that in the case where \(K \leq 0\), if
    \[ \inf\{\lla R_N(p)(X,Y)Y, X\rra: p\in \Sigma^2\ \text{and }X,Y \in T_p\Sigma\ \text{is orthonormal basis}\} < 0, \]
    where \(R_N\) represents the Riemann curvature tensor of \(N^n\), then, based on the Gauss equation, \(\Sigma^2\) must be a cylinder over a geodesic. In particular, if \(N^n\) has negative sectional curvature, then any minimal surface has \(K \leq 0\) and if \(|T|\) is not close to zero then \(\Sigma^2\) must be a cylinder over a geodesic (as in Corollary \ref{cor:trivialCor}). 

    In the work \cite{EspinarRosenbergBernstein}, results related to surfaces \(\Sigma^2\) with constant mean curvature in a product space \(N^2\times \Euc\), where \(N^2\) is a complete Riemannian manifold with non-constant curvature, were established. The authors assumed that the angle function retains a fixed sign on \(\Sigma^2\) and then categorized these surfaces based on the infimum of the Gaussian curvature of $N^2$ along the points of the \(\Sigma^2\) onto \(N^2\).

    \section{Final remarks}\label{sec:Remarks}
    
    In \cite{EspinarRosenbergCMCHomogeneous}, they concluded that beyond the surface being flat, if the norm of the Codazzi tensor is constant, then the angle will also be constant (as shown in Theorem 2.3 of \cite{EspinarRosenbergCMCHomogeneous}). This is in contrast to the higher codimension case. In comparison to \cite{AlencarTribuzyHopfMxR}, their conclusion is that the surface must necessarily be flat. They also posed the problem of classifying such surfaces (as mentioned in Remark 3 of \cite{AlencarTribuzyHopfMxR}). This classification was addressed later in \cite{HouQiuFlatProductConstantAngle}, under the assumption that the angle is constant. However, without making the assumption of a constant angle, achieving a classification is not expected without some control over the extrinsic curvature. An alternative way to see this is through the equation $$\frac{1}{2}\Delta|T|^2 = |A_{\eta}|^2 + \kappa |T|^2 (1 - |T|^2) - \sum_{i>1} |A_i T|^2,$$ where, for hypersurfaces, the last term vanishes. Under the hypotheses of Theorem \ref{thm:PMCSurfacesProducts}, it is a straightforward to deduce that $|T|^2$ is a subharmonic function on a parabolic surface, implying its constancy. With respect to the codimension, while both Theorems \ref{thm:PMCSurfacesProductsMu} and \ref{thm:PMCSurfacesProducts} are stated for arbitrary codimensions, \cite{AlencarTribuzyHopfMxR} offers specific results for codimension reduction that can be applied in this scenario.

    It is also worth noting that the items in Theorem \ref{thm:MainCodazzi} become equivalent when the operator $S$ is nonsingular. This is because assuming that $S$ is a Codazzi operator, defining the new metric on the surface $\lla \cdot, \cdot \rra = \lla S \cdot, S\cdot \rra$ for $X,Y \in T\Sigma$ ensures that $S^{-1}$ becomes a Codazzi operator in relation to the new connection. This naturally inverts the statement mentioned in Theorem \ref{thm:MainCodazzi}. This observation can also be directly applied with Simons' formula, as given by equation \eqref{eq:LnS}, to obtain the expression $\Tilde{\Delta} \ln |S| = \frac{4K}{|S|^2}$, which we utilized in Theorem \ref{thm:MainCodazzi}.

    The results in this manuscript are presented for orientable surfaces, but one can always work with the orientable double cover.

\section{Acknowledgements}
The author was supported by the Research Foundation-Flanders (FWO) and the Fonds de
la Recherche Scientifique (FNRS) under EOS Project G0H4518N.

\bibliographystyle{amsalpha}
\bibliography{bibliography}

\end{document}